\def\year{2020}\relax
\documentclass[letterpaper]{article} % DO NOT CHANGE THIS
\usepackage{aaai20}  % DO NOT CHANGE THIS
\usepackage{times}  % DO NOT CHANGE THIS
\usepackage{helvet} % DO NOT CHANGE THIS
\usepackage{courier}  % DO NOT CHANGE THIS
\usepackage[hyphens]{url}  % DO NOT CHANGE THIS
\usepackage{graphicx} % DO NOT CHANGE THIS
\usepackage{graphicx}  % DO NOT CHANGE THIS

\usepackage[utf8]{inputenc}
\usepackage{mathptmx} % assumes new font selection scheme installed
\usepackage{amsmath,amssymb} % assumes amsmath package installed
\usepackage{dsfont}
\usepackage{verbatim}
\usepackage{multicol,longtable}
\usepackage{epsfig}
\usepackage{txfonts}
\usepackage{calc}
\usepackage{subcaption}
\usepackage{enumerate}
\usepackage{bm}
\usepackage{color}
\usepackage{url}
\usepackage{algorithm}
\usepackage[noend]{algpseudocode}
\usepackage[table]{xcolor}

\title{AI Enhanced Control Engineering Methods \thanks{This paper was presented at the 2020 \textbf{Physics-guided AI to Accelerate Scientific Discovery (PGAI-AAAI-20)} Workshop.}}
\author{Ion Matei, Raj Minhas, Johan de Kleer and Alexander Felman\\ % All authors must be in the same font size and format. Use \Large and \textbf to achieve this result when breaking a line
Palo Alt Research Center\\ %If you have multiple authors and multiple affiliations
3333 Coyote Hill Road\\
Palo Alto, California 94304\\
imatei@parc.com, rajminhas@gmail.com, dekleer@parc.com, afeldman@parc.com % email address must be in roman text type, not monospace or sans serif
}

\begin{document}

\maketitle

\begin{abstract}
AI and machine learning based approaches are becoming ubiquitous in almost all engineering fields. Control engineering cannot escape this trend. In this paper, we explore how AI tools can be useful in control applications. The core tool we focus on is automatic differentiation. Two immediate applications are linearization of system dynamics for local stability analysis or for state estimation using Kalman filters. We also explore other usages such as conversion of differential algebraic equations to ordinary differential equations for control design. In addition, we explore the use of machine learning models for global parameterizations of state vectors and control inputs in model predictive control applications. For each considered use case, we give examples and results.
\end{abstract}

\section{Introduction}
Two modern pillars of AI, machine learning (ML) and deep learning (DL), with their models, platforms and training algorithms, are making their way into almost every engineering field. Rather than trying to find reasons why such tools and platforms cannot beat existing control methods, we explore how tools commonly used in ML/DL can be of help in control applications. Although we will refer to some DL models, the models used in our examples are minuscule compared to the deep learning models with millions of parameters that can be found in image classification applications \cite{5537907}, for example. Instead, we focus on one of the main tools provided by DL platforms, namely  \textit{automatic differentiation} (AD). It should not come as a surprise that there are quite a few immediate control applications. We can use such a tool to effortlessly compute Jacobian matrices of vector fields, providing easy to generate local linearizations of system dynamics. This use case is particularly interesting when dealing with large scale system models that are expressed in some modeling language (e.g., Matlab/Simulink or Modelica). Combined with a model transformation module that generates a representation compatible with the programming language on which the ML/DL platform resides (e.g., Python), AD can be a powerful tool for local stability analysis and control design. System dynamics Jacobian matrices are also used in state estimation algorithms. The extended Kalman filter \cite{Kalman,4501892} makes use of such matrices to track the states of nonlinear systems. Although there are other alternative approaches for state estimation (e.g., particle filter \cite{Arulampalam02atutorial}), the simplicity and explainability of Kalman  continue to make them good candidates for state estimation. We will also explore other less obvious usages of ML/DL platforms. Through an example, we will demonstrate how we can transform differential algebraic equations (DAEs) into ordinary differential equations (ODEs). The motivation is rather simple. We would like to extend the class of control design and state estimation approaches to systems whose dynamics are expressed as DAEs (e.g., chemical processes or mechanical systems). Structural analysis of DAEs combined with learning platforms is at the core of the DAE to ODE transformation. Another application of AD is model predictive control (MPC) \cite{GARCIA1989335}. We show that global parameterizations of state and control variables combined with automatically computed state time derivatives provide an alternative to traditional local parameterizations based on collocation methods \cite{doi:10.1137/16M1062569}. In addition, automatically computed gradients of objective functions and Jacobian matrices of constrained functions can considerably ease the process of formulating and solving the constrained nonlinear program behind the MPC formulation.

We do not favor any particular ML/DL platform. We do note though that platforms that integrate seamlessly with scientific programming languages such as Python considerably reduce the learning curve. Here is summary of a few characteristics of ML/DL platforms that we consider to be relevant in control applications:
\begin{itemize}
  \item \textit{Reduced learning curve}:  Users would rather prefer to leverage existing knowledge instead of having to learn fundamentally different frameworks or paradigms. New platforms that use non-standard, custom designed interfaces and constructs typically require a significant effort to get a handle on their usage. As a consequence, their chances of being widely adopted are small.
  \item \textit{Integration with scientific computing packages}: Formulation of optimization problems can include a variety of operators ranging from common mathematical operations (e.g., multiplication, addition) to more complex ones (e.g., matrix inversion, Fourier transforms). To effectively use AD the derivative operators must be wrapped around the aforementioned type of operators.
  \item \textit{The ability to generate compiled functions}: Time per iteration determines the feasibility of real time implementation of optimization algorithms. A way to achieve feasibility is the implementation of optimization algorithm in  programming languages that can generate compiled machine code (e.g., C/C++). If we use AD to generate functions that evaluate gradient/Jacobian, we need the ability to compile them and integrate them with the optimization algorithm compiled in machine code.
  \item \textit{Batch execution}: One major advantage of ML/DL models is their ability to evaluate and propagate multiple sample inputs. This is the case for neural networks which can efficiently process large batches of inputs. This is not enough when considering AD though. We need this property to be valid for functions representing gradients or Jacobian matrices as well. For example, in an MPC formulation we would like to be able to simultaneously process all time samples fed to the time derivative of the state variable. The alternative would be to use a ``for'' loop which is notoriously slow.
\end{itemize}
\textbf{Paper structure}: We first show a basic use of AD for local stability analysis and linear control design, and discuss an immediate extension to nonlinear state estimation. Next, we demonstrate how structural analysis of system dynamics can be coupled with regression problems to convert DAEs into ODEs, thereby extending the class of stability and control design approaches that can be applied to DAEs. Finally, we show an additional use case for AD when applied to MPC. In particular, we use AD to for batch evaluation of state time derivatives and for generating gradients/Jacobian matrices of objective and constraint functions part of the nonlinear program of the MPC formulation.

To execute AD we, use  the Jax \cite{jax2018github} Python library. Jax can automatically differentiate native Python and NumPy functions. It can compile and run NumPy programs on GPUs and TPUs; use just-in-time (JIT) compilation of Python functions for higher execution efficiency; and supports automatic vectorization. The latter is particularly useful since it allows batch evaluations of gradient functions. In addition, we use the {\tt scikit-learn} Python library to solve the linear regression problem formulated in the DAE to ODE conversion process.

\section{Automatic differentiation enabled stability analysis and control design}
\label{lab: stability}
Arguably, the  most basic approach for stability, controllability and observability analysis of nonlinear systems is linearization around an equilibrium point followed by computation of eigenvalues  and Grammian matrices. Linearized versions of dynamical systems are used to design linear controllers that enable local stability around equilibrium points. For small scale dynamical models (e.g., ODEs), the linearization can be done manually or through the use of symbolic calculus. For large scale systems though, symbolic calculus is not applicable. Hence we are left with the numerical approximations that typically induce numerical error. AD can be a solution when coupled with a model transformation module. The transformation module is necessary since dynamical models are expressed in different modeling languages. In particular, we are interested in models expressed in physics-based modeling languages (e.g., Modelica \cite{Fritzson15}) since they embed physical explainability. Thus, the transformation module is necessary for two reasons: (i) to perform simplifications of the physics-based model, and (ii) to generate a representation of the dynamical model compatible with the language where automatic differentiation resides (e.g., Python).

In what follows, we describe a workflow we used to transform physics-based models expressed in Modelica into Python code on which automatic differentiation can be applied. We use Jax \cite{jax2018github} for AD since it is compatible with  Python code expressed using Numpy and Scipy packages. We do not claim that our approach is unique, being determined by the authors' expertise and experience with a set of technologies related to physics-based modeling languages.
As a working example, we use the inverted pendulum dynamical model since it is both nonlinear and simple enough to show results at different stages of the transformation process. Note that the process can be readily generalized to much larger models. The dynamical model of the inverted pendulum is given by the following ODE:
{\small \begin{eqnarray}
\label{eq:04291759}
\ddot{x} &=& \frac{(I + ml^2)(b\dot{x}-F-ml\dot{{\theta}}^2\sin{{\theta}})-m^2l^2\frac{g}{2}\sin(2{\theta})}{m^2l^2\cos{{\theta}}^2-(m+M)(I + ml^2)}\\
\label{eq:04291760}
\ddot{\theta} & = & \frac{ml(F\cos{{\theta}}+ml\dot{\theta}^2\frac{g}{2}\sin{2{\theta}}-b\dot{x}\cos{{\theta}}+(m+M)g\sin{{\theta}})}{m^2l^2\cos{{\theta}}^2-(m+M)(I + ml^2)},
\end{eqnarray}}%
where $I$ is the pendulum moment of inertia, $M$ is the mass of the cart, $m$ is the pendulum mass, $b$ is the friction coefficient of the cart, $l$ is the length of the pendulum, $F$ is the force applied to the cart, $\theta$ is the pendulum angle and $x$ is the cart position. Let $z  =[x,\dot{x},\theta, \dot{\theta}]$, $u=F$ and $f(z,u)$ denote  the state vector, the control input, and system vector field (i.e., right hand side of  (\ref{eq:04291759})-(\ref{eq:04291760})), respectively.

\begin{figure}[htp!]
  \centering
  \includegraphics[width=0.98\linewidth]{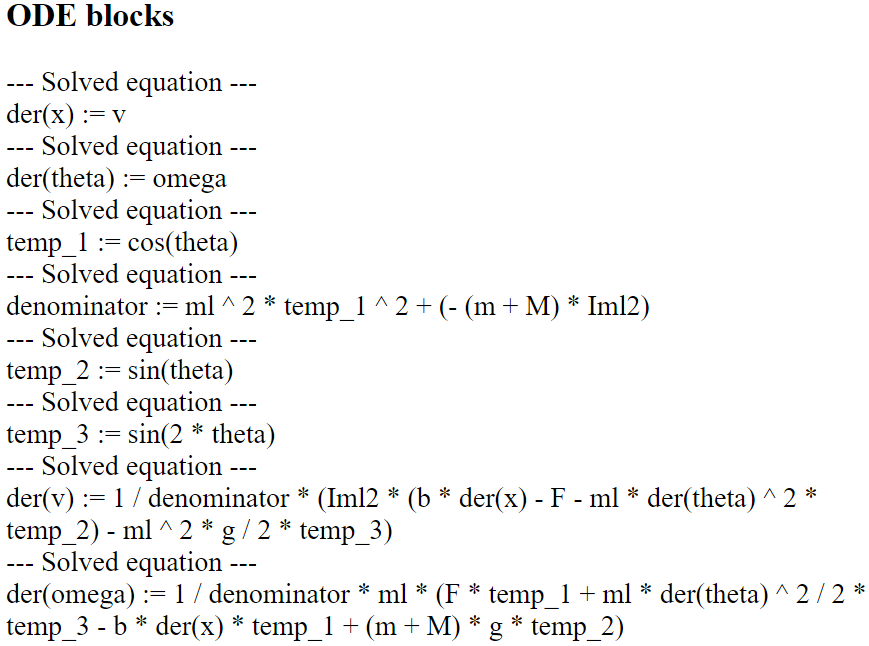}
  \caption{BLT form for the inverted pendulum (\ref{eq:04291759})-(\ref{eq:04291760}). The causal equations for computing the state derivatives are depicted.}
  \label{fig:BLT form inv pend}
\end{figure}

We used a component based approach to construct an inverted pendulum model in Modelica. Originally, the model is represented as a DAE. The DAE in this case was simple enough to be manually converted into an ODE by symbolically solving a linear system  of equations  with the state variables representing the unknowns. The model transformation process has the following steps:%

\textit{Modelica model analysis}:  The model is structurally analyzed  and a matching problem is solved that assigns variables to equations \cite{Casella,Casella11}. Since we are assuming that the model is expressed as an ODE, we do not have to worry about index reduction algorithms, as is the case with DAEs. The simplified equations are re-arranged to make their solution as simple as possible: the system is represented in a Block-Lower Triangular (BLT) form.
The BLT transformation problem is solved using Tarjan's algorithm. An additional step is applied: when possible, linear equations are solved symbolically. In the case of DAEs, additional processing steps are applied as discussed in a later section. We use the JModelica \cite{jmodelica} compiler to process the Modelica model and generate a simplified version of it which depicts a system of equations for computing the state derivatives. In particular, when compiling the Modelica model we enable the ``diagnostic" feature that generates the compilation results (i.e., the BLT form) as HTML files. A snapshot of the HTML file showing the BLT form for the inverted pendulum is shown in Figure \ref{fig:BLT form inv pend}. A set of assignment operations for computing the state derivatives are depicted, where the variable names are self-explanatory.

\textit{BLT form parsing}: We process the BLT form using an HTML parser. In particular we use the {\tt BeautifulSoup}  Python library to extract the BLT from the HTML file generated by JModelica. We extract the type of the quantities present in the model to create the state vector and list of parameters. Next we analyze the equations, extract the state variables, and convert the system of equations from strings to symbolic representations using the {\tt SymPy} Python library. Finally, we convert the symbolic representation into a representation compatible with AD. In particular, we generate a function that takes as arguments the state vector, the input, and system parameters and returns the time derivative of the state vector. In other words, this function evaluates the vector field at a particular state and input. %The Python code corresponding to the vector field function can be seen in Appendix \ref{app 2}. 
It is this function that is used by Jax for AD. %Note that the state vector is generated by taking individual state derivatives rather than index-based assignment since automatic differentiation does not support the latter.

The linearization of the nonlinear system dynamics is immediate. The derivative of the vector field function with respect to the state is obtained by calling the forward Jacobian Jax function, i.e., {\tt dfdz = jax.jacfwd(f)}. In a similar way we compute the derivative with respect to the input. To generate the matrices $A$ and $B$ of the linear system $\dot{z} = Az+Bu$, all we need to do is to evaluate the two Jacobian functions at the equilibrium points $x_{ss}$ and $u_{ss}$, that is {\tt A=dfdz(x\_{ss}, u\_{ss},params)} and {\tt B=dfdu(u\_{ss}, x\_{ss},params)}. Using the dictionary of parameters {\tt params=\{'M':0.5,'m':0.2,'b':0.1,'I':0.006,\\'g':9.81,'l':0.3,'ml':0.06,'Iml2':0.024\}}, we readily obtain the matrices
$$A=\left(
\begin{array}{cccc}
0 & 1 & 0 & 0\\
0 & -0.1818 & 2.6727 & 0\\
0 & 0 & 1 & 0\\
0 & -0.4545 & 31.18 & 0
\end{array}
\right)\ B^T=[0, 0.1818, 0, 4.5454]$$
These matrices can now be used for controllability analysis and for linear  control design.

\textbf{Remark:} We showed above how AD can be used for control design. An immediate extension of the use case is state estimation for nonlinear systems. Indeed, the extended Kalman filter \cite{Kalman,4501892} can definitely benefit from the Jacobian matrices generated through automatic differentiation, which can be used to linearize the sensing model as well. Instead of calling the Jacobian matrices at equilibrium points, we evaluate them at current state estimates, that is {\tt A=dfdz(x\_{hat}, u)}, where {\tt x\_{hat}} and {\tt u} denote the current state estimate vector and the control input.

\section{Conversion of DAEs to ODEs for control design}
\label{lab:conversion}
In this section we describe a method for converting DAEs to ODEs for the purpose of stability analysis and control design. We first talk about the method's steps, followed by a control design example.
\subsection{Method description}
A large class of mechanical, fluid or thermal physical systems cannot be represented as ODEs, but as DAEs of the form:
\begin{equation}
\label{eq:04051430}
F(\dot{x},x,u)=0.
\end{equation}
This is an implicit representation, where $\dot{x}$ cannot be represented analytically as a function of the state and input. Such systems can sometimes be brought into an ODE form by using index reduction techniques. When index reduction fails to transform a DAE into an ODE, we can again use Newton-Rapshon algorithm to compute the solutions of algebraic loops, or nonlinear equations that cannot be solved symbolically. This operation is numerically expensive since it requires the inversion of a Jacobian matrix which has $O(n^3)$ complexity, where $n$ is the number equations.
Arguably, designing controllers and state estimators for systems represented as DAEs is more challenging \cite{Daoutidis2015,KUMAR1997393} and typically aimed at particular types of classes of DAEs (e.g., linear form). We would like to augment the class of control and estimation methods available to DAEs by converting them into ODEs. The conversion process has two steps: (i) index reduction based on structural analysis, (ii) replacing the algebraic loops involving state derivatives with \emph{surrogate models} explicitly computing their solutions.
In what follows, we demonstrate the conversion process through an example from the chemical engineering domain. We consider a reactor with fast and slow reactions \cite{KUMAR1997393}, and our objective is to transform its DAE model into an ODE in order to be able to apply the ODE control theory to DAE systems as well. The considered DAE model describes two first order reactions, where a reactant $A$ is fed at a volumetric flowrate $F_i$ and concentration $C_{Ai}$, and the reactions $A\rightleftharpoons B$, and $B\rightarrow C$ occur in series. The reversible reaction $A\rightleftharpoons B$ is much faster than the irreversible one, and therefore is practically at equilibrium. The DAE governing the behavior of the two reactions is given by:
{\small\begin{eqnarray}
\label{eq:04111837}
\frac{dV}{dT} &=& F_i-F,\\
\label{eq:04111838}
\frac{dC_A}{dt} &=& \frac{F_i}{V}(C_{Ai}-C_A)-R_A,\\
\label{eq:04111839}
\frac{dC_B}{dt} &=& -\frac{F_i}{V}C_B+R_A-R_B,\\
\label{eq:04111840}
\frac{dC_C}{dt} &=& -\frac{F_i}{V}C_C+R_B,\\
\label{eq:04111841_}
0 &=& C_A-\frac{C_B}{K_{eq}},\\
\label{eq:04111841}
0 &=& R_B-K_BC_B,
\end{eqnarray}}
where $V$ is the volume of liquid in the reactor, $C_A$, $C_B$, $C_C$ are the molar concentrations of the corresponding components, and $R_A$ and $R_B$ are the reaction rates for the reversible and irreversible reactions, respectively. The control inputs for the reactions are the volumetric flowrates $F_i$ and $F$. The reactor model is a DAE with index two. The control objective is to maintain the liquid volume $V$ and the concentration of product $B$ ($C_B$) and some desired reference values.

Before solving a DAE, the solvers extract the structural information of the model and a matching problem is solved that matches variables with equations \cite{Casella,Casella11}. If matching fails (i.e., not all variables are associated to equations), index reduction algorithms, such as Pantelides's algorithm, are used to create new equations through repeated differentiation. Next, the equations are re-arranged to make their solution as simple as possible: the system is represented in a Block-Lower Triangular (BLT) form. The solution of the complete system is then reduced to the solution of many smaller sub-systems, whose dimensions correspond to that of the blocks on the main diagonal. The BLT transformation problem is solved using Tarjan's algorithm. Two additional steps are applied: when possible, linear equations are solved symbolically, and non-linear equations are simplified by the \emph{tearing} method. In the tearing step, part of the variable in a block that corresponds to a set of nonlinear equations, is solved in terms of the remaining variables. As a result, the Newton-Raphson algorithm is applied to a smaller number of variables.

\begin{figure}[htp!]
  \centering
  \includegraphics[width=0.8\linewidth]{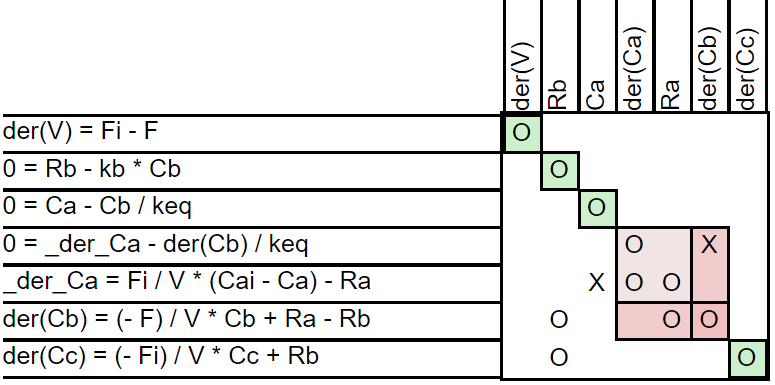}
  \caption{BLT form for the reactor model (\ref{eq:04111837})-(\ref{eq:04111841}).}
  \label{fig:BLT form}
\end{figure}

The structural analysis applied to the DAE (\ref{eq:04111837})-(\ref{eq:04111841}) generates the BLT form shown in Figure \ref{fig:BLT form}.
The index reduction algorithm introduced the ``dummy'' variable $\dot{C}_A$ that decouples the computation of $C_A$ from its time derivative $\frac{dC_A}{dt}$. The unknown variables of the diagonal block are $\dot{C}_A$, $\frac{dC_B}{dt}$ and $R_A$ - the solution is computed as a function of the system variables $F_i$, $F$, $V$, $C_A$, $C_B$ and $R_B$. A quick exploration of the BLT form shows us that we can eliminate variable $R_B$ since it depends linearly on $C_B$. Let $y^T=\left[\dot{C}_A,\frac{dC_B}{dt},R_A\right]$ and $u^T=\left[F_i, F, V, C_A, C_B\right]$. Our aim is to determine an explicit model $y = f(u;w)$ that computes the solution of the diagonal block, where $w$ is a vector of parameters.  In this example, we use a more traditional ML approach, where we engineer a set of features and formulate a regression problem. Alternatively, we can use a NN to model $f(u;w)$.  We generate the training data by randomly sampling the input space and solving the BLT form diagonal block, expressed as:
{\small\begin{eqnarray}
\label{eq:04111904}
0 &=& \dot{C}_A - \frac{1}{K_{eq}}\frac{dC_B}{dt},\\
\label{eq:04111905}
\dot{C}_A &=& \frac{F_i}{V}(C_{Ai}-C_A)-R_A,\\
\label{eq:04111906}
\frac{dC_B}{dt} &=& -\frac{F}{V}C_B+R_A-K_{B}C_B
\end{eqnarray}}%
Note that (\ref{eq:04111904})-(\ref{eq:04111906}) represents the ground truth and we can generate arbitrarily large training data sets. When sampling the input space we consider the underlying physical constraints that apply to the system variables: non-negative volumetric flowrates, non-negative liquid volume, non-negative component concentrations. We define the following feature vector $x=[$$x_1$,$x_2$,$x_3$,$x_4$, $x_1x_4$, $x_2x_3$,$x_1x_2$, $x_1x_3$,$x_2x_4$, $x_3x_4]^T$, where $x_1=F_i/V$, $x_2=F/V$, $x_3=C_B$ and $x_4=C_A$. We consider the linear regression model $y_j = w_j^T x$, where $y_1 = C_A$ and $y_2 = C_B$ and $y_3 = R_A$. Table \ref{tab:04111919} shows the results of training the parameters of the regression model.
 \begin{figure}[htp!]
  \centering
  \includegraphics[width=0.9\linewidth]{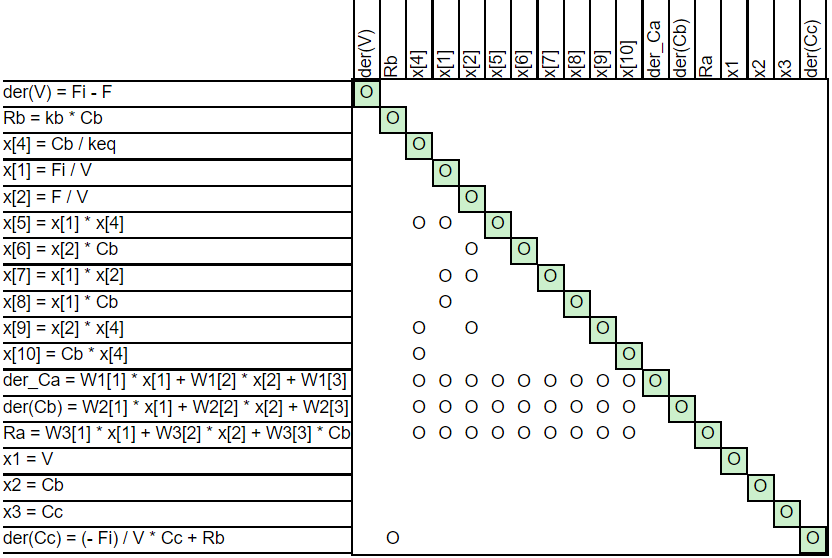}
  \caption{BLT form for the reactor model with explicit solution of the diagonal block.}
  \label{fig:BLT form ode}
\end{figure}

\textbf{Remark:} For this particular example, the system of equations (\ref{eq:04111904})-(\ref{eq:04111906}) is linear. For small to medium sized linear systems, we can use symbolic calculus to compute their solutions. For large linear or nonlinear systems of equations, symbolic calculus may fail, and we are forced to numerically learn a model for their solutions, as we did in our example.
\begin{table*}[htp!]
\caption{Parameters of the regression model}
\centering
\begin{tabular}{ |c|c| }
 \hline
$w_1^T$ & [3.333e+01,  5.432e-14, -2.000e-01,  6.418e-14,
  -6.666e-01, -6.666e-01,  1.352e-15, -8.850e-16,
  -1.612e-16,  1.123e-15]  \\
   \hline
 $w_2^T$ & [1.666e+01,  2.7164e-14, -1.000e-01,  3.209e-14,
  -3.333e-01, -3.333e-01,  6.761e-16, -4.425e-16,
  -8.061e-17,  5.619e-16] \\
   \hline
$w_3^T$ & [1.666e+01, -2.377e-14,  2.000e-01, -1.359e-14,
  -3.333e-01,  6.666e-01, -7.721e-17, -1.704e-16,
  -1.270e-15,  1.406e-15] \\
 \hline
\end{tabular}
\label{tab:04111919}
\end{table*}
We replace the equations in the block diagonal with the regression model that computes the solution of the block, and again generate BLT form of the system model, as shown in Figure \ref{fig:BLT form ode}.

We note that the diagonal block has disappeared and the BLT form is a pure triangular form. More importantly, the original DAE model can now be expressed as an ODE with three state variables:
\begin{eqnarray}
\label{eq:04111943}
\frac{dV}{dt} &:=& F_i-F,\\
\label{eq:04111944}
C_A &:=& \frac{C_B}{K_{eq}}\\
\nonumber
x &:=& \left[\frac{F_i}{V},\frac{F}{V},C_B,C_A, \frac{F_i}{V}C_A, \frac{F}{V}C_B,\frac{F_iF}{V^2}, \right.\\
\label{eq:04111945}
& & \left.\frac{F_i}{V}C_B,\frac{F}{V}C_A, C_AC_B\right]^T,\\
\label{eq:04111946}
\frac{dC_B}{dt} &=& w_2^Tx,\\
\label{eq:04111947}
\frac{dC_C}{dt} &=& \frac{F_i}{V}C_C+K_BC_{B}.
\end{eqnarray}
By substituting the variable $C_A$ in $x$ as a function of $C_B$, we obtain the familiar form of an ODE.

\subsection{Control design}
\begin{figure}[htp!]
  \centering
  \includegraphics[width=0.9\linewidth]{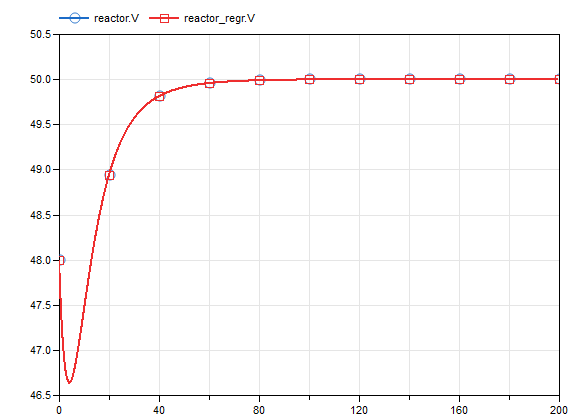}
  \caption{Evolution of the variable $V$ under feedback control, with initial states $V(0) = 48$ liters, $C_B(0)=12$ mols/liter, and $C_C(0) = 19$ mols/liter.}
  \label{fig:V closed loop}
\end{figure}

\begin{figure}[htp!]
  \centering
  \includegraphics[width=0.9\linewidth]{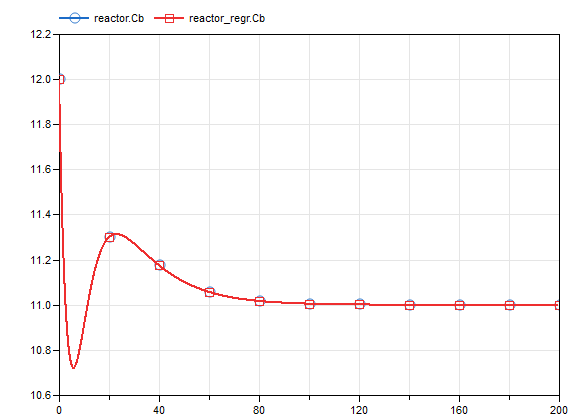}
  \caption{Evolution of the variable $C_B$ under feedback control, with initial states $V(0) = 48$ liters, $C_B(0)=12$ mols/liter, and $C_C(0) = 19$ mols/liter.}
  \label{fig:Cb closed loop}
\end{figure}

\begin{figure}[htp!]
  \centering
  \includegraphics[width=0.9\linewidth]{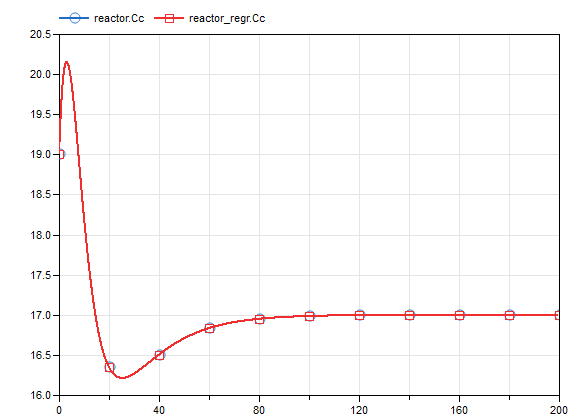}
  \caption{Evolution of the variable $C_C$ under feedback control, with initial states $V(0) = 48$ liters, $C_B(0)=12$ mols/liter, and $C_C(0) = 19$ mols/liter.}
  \label{fig:Cc closed loop}
\end{figure}

The ODE form enables linearization of system dynamics at an equilibrium point that, in turn, enables controllability and observability analysis. We first determine the equilibrium point. The control objective is to stabilize the system at the equilibrium point $C_B^{ss} = 11 $ mol/liters and $V^{ss}=50$ liters. Using (\ref{eq:04111943})-(\ref{eq:04111946}),  we can solve for the remaining steady state variables and inputs: $F_i^{ss}=9.705$ liters/sec, $F^{ss}=9.705$ liters/sec, $C_C^{ss}=17.0$ mols/liter. There are a plethora of control techniques that can be used to stabilize the system \cite{khalil2002nonlinear}. As an example, we design a linear controller  for the linearized system model at the equilibrium point. We use the AD feature of Jax \cite{jax2018github} to generate the linear model. In particular, given the nonlinear ODE $\dot{x} = f(x,u)$, we compute the Jacobian matrices $A = \frac{\partial f}{\partial x}(x^{ss}, u^{ss})$ and $B = \frac{\partial f}{\partial u}(x^{ss}, u^{ss})$. The linearized system is given by $\frac{d \Delta x}{dt}=A \Delta x + B\Delta u$, where $\Delta x = x-x^{ss}$ and $\Delta u = u-u^{ss}$, with
$$A = \left[\begin{array}{ccc}
0 & 0 & 0\\
-0.022 & -0.2941 & 0\\
0.0660 & 0.3 & -0.1941
\end{array}
\right],$$
$$B = \left[\begin{array}{cc}
1 & -1 \\
0.1867 & -0.0733\\
-0.3400 & 0
\end{array}
\right],
$$
where the system parameters were chosen as: $K_{eq} = 0.5$, $K_B = 0.3$ and $C_{Ai} = 50$ mols/liter.
We use a pole placement approach, where we generate the matrix $K=[k_1,k_2]^T$ such the eigenvalues of the matrix $A-BK$ are at $p=[-0.1, -0.5, -0.3]$. The control inputs $F_i$ and $F$ are computed as $F_i = k_1^T(x-x^{ss})+F_i^{ss}$ and $F = k_2^T(x-x^{ss})+F^{ss}$, where $x = [V,C_B, C_C]^T$. Figures \ref{fig:V closed loop}, \ref{fig:Cb closed loop} and \ref{fig:Cc closed loop} show the evolution of the state variables under feedback control, as they stabilize at the desired equilibrium point. Each figure shows the state for both the linearized model and for the DAE model when the inputs are generated by the linear controller.

We tested the controller for different initial conditions, empirically observing that the closed loop system has a large attractor. In general, we can expect to introduce uncertainty in the system model by replacing diagonal blocks with causal models. %Such uncertainty can be quantified, as described in the case of the pendulum system.
Once the uncertainty is quantified, we can use robust control techniques \cite{morari1989robust} to account for it and to ensure the stability of the closed loop system.

\section{Global parameterization of system states and controls for MPC}
\label{lab:mpc}
In this section we show how we can use global parmeterization of state and input variables to compute time dependent control inputs over some finite time horizon. We will use the inverted pendulum system as an example. The nonlinear dynamics of the inverted pendulum were introduced  in (\ref{eq:04291759})-(\ref{eq:04291760})).

 MPC \cite{GARCIA1989335} is the standard approach for control of nonlinear systems and requires solving a nonlinear program with constraints. The optimization variables are the state trajectories and control inputs over time. The most important constraint is the system dynamics constraint that makes sure that the state trajectories respect the system dynamics. Other constraints can include setting the initial and final values of the state variables, or limiting the magnitude of the states and control inputs. 
 
 The MPC formulation assumes a discretization of the time domain $\mathcal{T}=\{t_0,t_1,\ldots, t_N\}$ and the optimization variables are the state and control variables at discrete time instances, that is, $z(t_k)$ and $u(t_k)$, for $k\in \{0,1,\ldots,N\}$. The key step in the formulation of the system dynamics constraint is the representation of the time derivatives $\dot{z}(t_k)$. Approaches based on trapezoidal  collocation or Hermite-Simpson collocation \cite{doi:10.1137/16M1062569} are typically used to transform the continuous dynamics into a discrete set of equality constraints. Such methods locally approximate the state and input trajectories while imposing equality constraints at the collocation points. We use a different approach to solve the MPC problem: a global representation of the state and input trajectories along with automatic differentiation to compute time derivatives. In particular, we choose neural networks to represent the state and input over time, that is $z:\mathds{R}\rightarrow\mathds{R}^4$ and $u:\mathds{R}\rightarrow\mathds{R}$.  We use automatic differentiation to explicitly compute the time derivative of the state vector. This way, we can avoid the need for collocation-like methods.
 Let $w$ denote the combined parameters of the state and control input representations. The constrained nonlinear program for the MPC formulation is given by:
\begin{eqnarray}
\label{eqn:04221335}
\min_{w} &  &\frac{1}{2}\|w\|^2\\
\label{eqn:04221336}
\textmd{subject to: } & & \dot{z}(t_k;w) = f(z(t_k;w),u(t_k;w), \forall t_k\in \mathcal{T},\\
\label{eqn:04221337}
& & {z}(t_0;w) = {z}_0,\\
\label{eqn:04221338}
& & {z}(t_N) = z_{N},\\
\label{eqn:04221339}
& & |u(t_k;w)|\leq U_{max},  \forall t_k\in \mathcal{T},
\end{eqnarray}
where $z_0$, $z_N$ and $U_{max}$ are the initial and final condition, and the maximum control magnitude, respectively. We can reduce the number of equality constraints by replacing   (\ref{eqn:04221336}) with an equivalent constraint 
\begin{equation}
    \label{eq:04301732}
    \frac{1}{N}\sum_{k=0}^N\|\dot{z}(t_k;w) - f(z(t_k;w),u(t_k;w)\|^2 = 0.
\end{equation}
Note that due to the parameterization of the state and input, we are not dependent on a particular discretization scheme. In fact, during the optimization procedure, we can randomly select time instances from the time domain. In turn, this will reduce overfitting of the optimization solution to a particular time discretization scheme. The optimization solution parameterization has the same flavor as the PDE solution parameterizations studied in \cite{alaradi2018solving}. We prefer to use a constrained optimization formulation that strictly enforces the equality constraints to avoid the painful exercise of properly scaling the different components of the loss function in case we add the constraints (as regularizers) in the loss function. 

We choose neural network representations for the state and control input for practical reasons: their evaluation can be done efficiently on GPU/TPUs. In other words, we can do batch executions to evaluate  the state, state derivatives and input for all  time samples. When choosing the optimization algorithm, we consider the fact that we can again use AD to compute gradients and Jacobian matrices of the loss and constrain functions. One algorithm for nonlinear programming that proves to be efficient in practical implementations and makes heavy use of gradients and Jacobian matrices is the sequential quadratic programming (SQP) algorithm \cite{NoceWrig06}. This algorithm solves a sequence of quadratic, convex optimization problems with equality constraints. AD deals with the linearization difficulties that sometimes plague the SQP algorithm. The SQP algorithm also uses second order information (e.g., Hessian matrix) of the loss function that can add additional complexities for large scale problems. In our case, we avoid these potential complexities by having a simple L2 norm cost, which has a constant Hessian matrix. For evaluating gradients and Jacobian matrices we use the Jax \cite{jax2018github} Python library. 
Jax enables a computationally efficient evaluation of constraints and their Jacobian matrices through batch executions. At each (outer) iteration of the SQP algorithm, we evaluate the state vector $z(t)$ and the control input $u(t)$ for all time samples at the same time. This is typical for NN models. Next, we evaluate the function $f(z(t),u(t))$, which can be computed for all time samples using a single statement, without the need of a {\tt for} loop. In the case of the state derivative, we enable batch execution through ``mapping''. First we generate the vector valued state derivative by using the Jacobian operator provided by Jax. Namely, we define {\tt dzdt = jax.jacfwd(z)} which returns a function that takes one argument, i.e. time, and generates a vector whose entries are the time derivatives of each state variable. To enable batching we apply mapping, which can be done in Jax through the statement {\tt dzdt\_vec = vmap(dzdt)}. We can now pass a vector of time samples to {\tt dzdt\_vec} and generate a matrix whose rows are the state derivatives at a particular time instance. Additional efficiency can be achieved by using the {\tt jit} (just-in-time) decorator which produces optimized kernels through accelerated linear algebra (XLA). There are several recent attempts to implement SQP algorithms on GPU platforms \cite{8253023,101145} that enables an increased efficiency in solving large scale MPC problems.

\begin{figure}
\begin{subfigure}[b]{.4\textwidth}
  \centering
  \includegraphics[width=\textwidth]{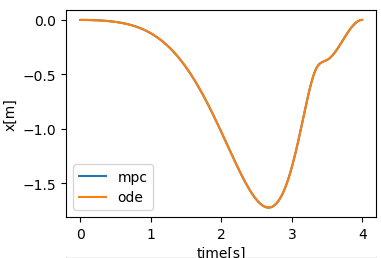}
  \caption{Cart position.}
  \label{fig:10c_potentials_t0}
\end{subfigure}%
\hfill
\begin{subfigure}[b]{.4\textwidth}
  \centering
  \includegraphics[width=\textwidth]{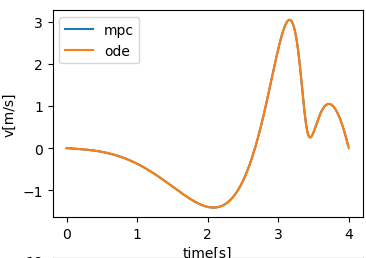}
  \caption{Cart velocity.}
  \label{fig:10c_potential_t1}
\end{subfigure}%
\hfill
\begin{subfigure}[b]{.4\textwidth}
  \centering
  \includegraphics[width=\textwidth]{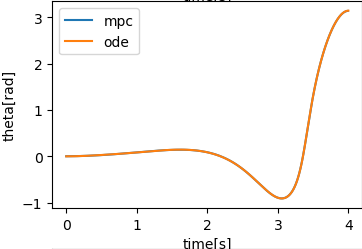}
  \caption{Pendulum angle.}
  \label{fig:10c_potential_t2}
\end{subfigure}%

\begin{subfigure}[b]{.4\textwidth}
  \centering
  \includegraphics[width=\textwidth]{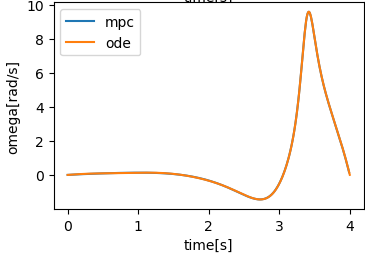}
  \caption{Pendulum angular velocity.}
  \label{fig:10c_potential_t3}
\end{subfigure}%
\caption{Inverted pendulum state variables: ``mpc'' label designates the solution of the MPC problem; ``ode'' designates the inverted pendulum ODE solution when using the control input generated by the MPC problem.}
\label{fig:variables}
\end{figure}

\begin{figure}[htp!]
  \centering
  \includegraphics[width=0.42\textwidth]{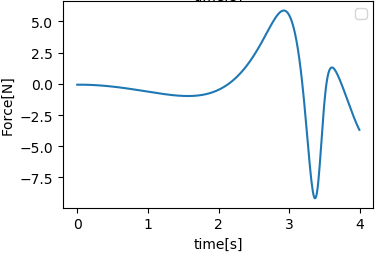}
  \caption{Force trajectory generated by the MPC problem}
  \label{fig:force}
\end{figure}

We integrated the gradients and Jacobian matrices of the loss and the constraints with the SQP algorithm from the  SciPy optimization package to learn a control input that takes the pendulum from an initial condition $z=[0,0,0,0]$ (pendulum pointing downwards) to a final condition $z = [0,0,\pi,0]$ (pendulum pointing upwards). In other words, we start far from the unstable equilibrium point. We limited the control magnitude to 10 N. We ran the algorithm for 2000 iterations with uniformly random initial values for $w$ from the interval $[0,10^{-2}]$. We used relatively small NNs, having one hidden layer of size 30 and {\tt tanh} as activation function.  As a result, the state and control input are  represented using 184 and 91 parameters, respectively.  These modest NN proved to be sufficient to learn the solution of the nonlinear program. More importantly, these numbers are independent on the number of time samples.  Figures \ref{fig:variables} show the trajectories of the inverted pendulum state variables. The ``mpc'' label designates the solution of the MPC problem, while ``ode'' label designates the solution of inverted pendulum ODE solution when using the control input over time generated by the MPC problem, but not the NN parameterization of the state. Figure \ref{fig:force} shows the force that needs to be applied to the cart to bring the pendulum to the desired final position. The figures show the results when starting from the initial condition $z_0$. The optimization problem is repeated at each time step successively to account for disturbances. Hence, it is imperative to do the  real time implementation in languages that generate machine code. As a consequence, we need the AD process to generate machine code compiled gradients and Jacobian matrices for fast control computations.

\section{Conclusions}
In this paper we showed how modern AI tools can help control engineers. We presented three examples. In the first example automatic differentiation coupled with model transformation was used to linearized nonlinear dynamical systems. Such linearization can be used for stability analysis and design of controllers and state estimators. In the second example we showed how we can transform DAEs to ODEs. Structural analysis together with learning regression models was at the core of the transformation. In the last example we showed how we can globally parameterized state and input vectors for MPC formulations. Automatic differentiation enhanced with batch execution enable efficient evaluations of state time derivatives, loss and constraint functions. In addition, automatically computed gradients and Jacobian matrices of loss and constraint functions removes the need of numerical approximations and hence prevents error accumulations during the optimization process. Combined with GPU/TPU implementations of constrained optimization algorithms, this mix of methods and technologies has the potential to provide solutions for solving large scale MPC problem in real time.

\section{ Acknowledgments}
This material is based, in part, upon work supported by the Defense
Advanced Research Projects Agency (DARPA) under Agreement No.
HR00111990027 and HR00111890037.

\bibliographystyle{aaai}
 \bibliography{references}

\end{document}